\documentclass[11pt]{article}

\date{}

\usepackage{url}
\usepackage{cite}
\usepackage{plain}

\usepackage[english]{babel}
\usepackage[leqno]{amsmath}
\usepackage{amssymb}
\usepackage{verbatim}
\usepackage[mathscr]{eucal}

\usepackage{amstext}
\usepackage{makeidx}

\usepackage{hyperref}
\hypersetup{colorlinks,%
citecolor=red,%
linkcolor=blue,%
urlcolor=black}

\makeindex

\newtheorem{theorem}{Theorem} 

\newtheorem{proposition}{Proprieta'}
\newtheorem{definition}{Definizione}
\newtheorem{notation}{Nota}
\newtheorem{ex}{Esercizio} 
\newtheorem{esempio}{Esempio}

\newcommand{\beq}{\begin{equation}} 
\newcommand{\eeq}{\end{equation}}

\newcommand{\bex}{\begin{ex}} 
\newcommand{\eex}{\end{ex}} 

\newcommand{\bese}{\begin{esempio}} 
\newcommand{\eese}{\end{esempio}} 

\newcommand{\bpro}{\begin{proposition}} 
\newcommand{\epro}{\end{proposition}}

\newcommand{\bthe}{\begin{theorem}} 
\newcommand{\ethe}{\end{theorem}}

\newcommand{\bnote}{\begin{notation}} 
\newcommand{\enote}{\end{notation}}

\newcommand{\bdefi}{\begin{definition}} 
\newcommand{\edefi}{\end{definition}} 

\newcommand{\bc}{\begin{center}} 
\newcommand{\ec}{\end{center}}

\newcommand{\mail}[1]{\href{unina:#1}{\texttt{#1}}}

\usepackage{pdfsync}

\author{ Monica De Angelis \thanks{Univ. of Naples  "Federico II", Dip. Mat. Appl. "R.Caccioppoli", \newline
  80125, Naples, Italy.
\newline\mail{modeange@unina.it},}
}

\title{A priori estimates for solutions of FitzHugh-Rinzel   system  }

\begin{document}
\maketitle

\begin{abstract}
 The FitzHugh-Rinzel system is able to describe some biophysical phenomena, such as bursting oscillations,
  and the study of its solutions can help to better understand  several behaviours  of  the complex dynamics of biological systems. We express the solutions  by means of  an  integral equation  involving the fundamental solution $ H(x,t)  $  related to a non linear integro-differential  equation. Properties  of  $ H(x,t)  $ allow us to obtain a priori estimates for solutions  determined in the whole  space, showing both the influence of the initial data and the source term. 
\end{abstract}

\section{Introduction }
\label{intro}

Mathematical biophysics models, such as the FitzHugh Nagumo system (FHN), play an important role in studying the nervous system, as they can help describe biophysical phenomena that are relevant to neuronal excitability.

The FHN consists of two differential equations that model several engineering  applications and there exist many scientific results and an extensive bibliography  in regard. \cite{i,dr13,moca21,r1,m13,glrs,drw}.  
 However, it has been noted that only using a reset or  adding noise, it is possible to evaluate  bursting phenomena. 
 This  phenomenon  occurs in a number of different cell types and  it consists of a  behaviour characterized by brief bursts of oscillatory activity alternating with periods of quiescence during which the membrane potential changes only slowly \cite{ks}.

Bursting phenomena  are becoming more and more important and their studies are increasing in many scientific fields (see, f.i.{\cite{2020} and references therein). For example, in  the  restoration of synaptic connections, it appears that some nanoscale memristor devices
 have the potential to reproduce the behavior of a biological synapse \cite{jnbbikem,clag}.
 This will lead in the future, especially  in  case of traumatic injuries, to the introduction
of electronic synapses to directly connect neurons.

A  model   that  seems to be more mathematically appropriate   for incorporating nerve cell bursting  phenomena is the FitzHugh Rinzel model (FHR). It derives from FitzHugh-Nagumo and, differently from FHN,  consists of three  equations  just to insert slow modulation of the current\cite{i,bbks,ws15,rt,r}. Indeed,  bursting oscillations  can be characterized by a system variable that periodically changes from an active phase of rapid spike oscillations to a silent phase.

  As for the FHR model, the following system is considered:

\begin{equation}
\label{11}
  \left \{
   \begin{array}{lll}
    \displaystyle{\frac{\partial \,u }{\partial \,t }} =\,  D \,\frac{\partial^2 \,u }{\partial \,x^2 }
     \,-\, w\,\,+y\,\, + f(u ) \,  \\
\\
\displaystyle{\frac{\partial \,w }{\partial \,t } }\, = \, \varepsilon (-\beta w +c +u)
\\
\\
\displaystyle{\frac{\partial \,y }{\partial \,t } }\, = \,\delta (-u +h -dy)
   \end{array}
  \right.
\end{equation}                                                     
 
\noindent where

 \begin{equation}                 \label{12}
f(u)= u\, (\, a-u \,) \, (\,u-1\,) \,\,\,\,\,\,\, \,\, (\,0\,<\,a\,< \,1\,),
\end{equation}

\vspace{3mm} \noindent and terms  $ \,\beta, \,\, c,\,\,d, \, \,h, \,\, \varepsilon, \,\,\delta \,\, $  are positive constants that  characterize the model's kinetics. The second order term with    $ D > 0 $   can be associated to the axial current in the axon, and  it derives from the Hodgkin-
Huxley theory for nerve membranes. Indeed, if $ b $ represents the axon diameter and
$ r_i $ is the resistivity, the spatial variation in the potential $ V $  gives the term $(b/4r_i )V_{xx}$
from which term $D u_{xx}$ derives (see f.i.\cite{m2}), and in \cite{2020}  an analysis on contribution due to this   term  has been developed. 
Furthermore, when
 the fast variable $ u$ simulates the membrane potential of a nerve cell, while the slow variable
 $ w$ and the super-slow variable $ y$ determine the corresponding number densities of ions, the model (\ref{11}) simulates the propagation of impulses from one neuron to another, and  studies  on solutions can help  in testing  the responses of the various models in neuroscience. 
 
 Several methods have been developed to find exact solutions related to  partial differential equations and     an extensive bibliography on the study of analytical behaviors exists (see,f.i \cite{dma18,Lg,df213,mda19,adc,mda13,k18}). The aim of this paper  is to  determine a priori estimates for the (FHR) solution  by means of suitable properties of  the fundamental solution $ H(x,t),  $ showing how the effects due to the initial perturbation are vanishing when $ t $ tends to infinity, and simultaneously, as time increases, the effect of the nonlinear source remains bounded.

The paper is organized as follows: in section 2 we define the mathematical problem and report some  of the results already proved in \cite{2020}, as well as other  results  well known in literature. In section 3,  some properties related to the  fundamental solution $ H(x,t) $ are obtained and, in a subsection, some relationships on convolutions which characterize the explicit solution, are highlighted. In section 4,  estimates on convolution  are proved and in section 5, the solution  is expressed   by means of these particular  convolution integrals. Finally, in section 6,  a priori estimates are showed.

\section{Mathematical considerations}

Indicating by $\,T\,$   an arbitrary  positive constant, let us  consider the set:
\begin{center}

$\  \Omega_T =\{(x,t) :  x \in \Re
, \  \ 0 < t \leq T \}.$
\end{center}

Moreover, if 

\begin{equation}      \label{ic}
u(x,0)\, =\,u_0 \,, \qquad w(x,0)\, =\,w_0  \qquad   y(x,0)= y_0,\qquad\qquad ( \,x\,  \in \Re\,)
\end{equation}

\noindent  represent the initial values, then  from  $(\ref{11})_{2,3}, $   one deduces:

\begin{equation}
\label{17}
\left \{
   \begin{array}{lll}
\displaystyle  w\, =\,w_0 \, e^{\,-\,\varepsilon \beta\,t\,} \,+\, \frac{c}{\beta}\,( 1- e^{- \, \varepsilon \, \beta \,t} )\,+ \varepsilon \int_0^t\, e^{\,-\,\varepsilon \,\beta\,(\,t-\tau\,)}\,u(x,\tau) \, d\tau 
  \\
  \\
\displaystyle  y\, =\,y_0 \, e^{\,-\,\delta\, d\,t\,} \,+\, \frac{h}{d}\,( 1- e^{ - \, \delta \,d  \,t} )\,- \delta \int_0^t\, e^{\,-\,\delta \,d\,(\,t-\tau\,)}\,u(x,\tau) \, d\tau.
 \end{array}
  \right.  
\end{equation}

\vspace{3mm} Besides, letting 

\begin{equation}      \label{14}
f(u)\, =\,-\,a\,u\, +\, \varphi(u) \quad with \quad  \varphi(u) \,=\, u^2\, (\,a+1\,-u\,)\quad  0<a<1, 
\end{equation}

 \vspace{3mm}\noindent  system (\ref{11}) becomes 
  
\begin{equation}
\label{15}
  \left \{
   \begin{array}{lll}
    \displaystyle{\frac{\partial \,u }{\partial \,t }} =\,  D \,\frac{\partial^2 \,u }{\partial \,x^2 } -au \,\,
     \,-\, w\,\,+y\,\, + \varphi(u ) \,  \\
\\
\displaystyle{\frac{\partial \,w }{\partial \,t } }\, = \, \varepsilon (-\beta w +c +u)
\\
\\
\displaystyle{\frac{\partial \,y }{\partial \,t } }\, = \,\delta (-u +h -dy),
   \end{array}
  \right.
\end{equation}

\vspace{3mm}\noindent and hence, when

\begin{equation}
\label{18}
\displaystyle F(x,t,u) =\varphi (u)  - w_0(x)  e^{- \varepsilon  \beta  t}+ y_0(x)  e^{- \delta  d  t}- \frac{c}{\beta}( 1- e^{- \varepsilon \beta t} )+\frac{h}{d}( 1- e^{ -\delta d t} )
\end{equation}

\vspace{3mm}\noindent denotes the source term, problem (\ref{15}) with initial conditions (\ref{ic}), can be modified into  the  following  initial value problem   $\,{\cal P}$:

	  \begin{equation}                                                     \label{19}
  \left \{
   \begin{array}{lll}
   \displaystyle  u_t - D  u_{xx} + au + \int^t_0 [ \varepsilon e^{- \varepsilon  \beta (t-\tau)}+ \delta e^{- \delta d (t-\tau)} ]u(x,\tau)  d\tau = F(x,t,u) \\    
\\ 
 \displaystyle \,u (x,0)\, = u_0(x)\,  \,\,\,\,\, x\, \in \, \Re. 
   \end{array}
  \right.
 \end{equation}

\vspace{3mm} In order to determine the solution of problem (\ref{19}), let us consider the following  functions:

   \begin{eqnarray} \label {22}
\nonumber \displaystyle  &H_1(x,t)\,=\, \, \, \frac{e^{- \frac{x^2}{4\,D\, t}\,}}{2 \sqrt{\pi  D t } }\,\,\, e^{-\,a\,t}\,
 \, 
 \\
 \\
&\nonumber \displaystyle - \frac{1}{2} \,\,    \,\,\,\int^t_0  \frac{e^{- \frac{x^2}{4 \, D\, y}\,- a\,y}}{\sqrt{t-y}} \,\, \, \frac{\,\sqrt{\varepsilon} \,\, e^{-\beta \varepsilon \,(\, t \,-\,y\,)}}{\sqrt{\pi \, D \,}} J_1 (\,2 \,\sqrt{\,\varepsilon \,y\,(t-y)\,}\,\,)\,\,\} dy,
   \end{eqnarray}

\begin{equation} \label{H2}
 \displaystyle H_2 =\int _0 ^t  H_1(x,y) \,\,e^{ -\delta d (t-y)} \sqrt{\frac{\delta y}{t-y}}   J_1( \,2 \,\sqrt{\,\delta \,y\,(t-y)\,}\,\,\, dy
 \end{equation}

\vspace{3mm}\noindent  
 where     $ J_1 (z) \,$  is the  Bessel function of first kind and order $\, 1.\,$

 \vspace{3mm}  In \cite{2020} it  has been verified  that function $ H(x,t): $ 

\begin{equation} \label{A16}
\displaystyle H = H_1  - H_2
\end{equation}

\noindent represents the fundamental solution of the parabolic 
 operator

 \begin{equation} \label{a1}
 {L }u\equiv  u_t - D  u_{xx} + au + \int^t_0 [ \varepsilon \,e^{-\, \varepsilon  \beta (t-\tau)}\,+ \delta e^{-\, \delta d (t-\tau)} ]u(x,\tau) \, d\tau, 
 \end{equation}

\vspace{3mm} \noindent and the following theorem has been proved:

\begin{theorem}  \label{th}

  In the half-plane $ \Re e  \,s > \,max(\,-\,a ,\,-\beta\varepsilon ,- \delta d\,)\,$    the Laplace integral  $\,{\cal L  }_t\, H \, \,$  converges absolutely for all  $\,x>0,\,$   and it results:
 \end{theorem}

\begin{equation}      \label{24}
\displaystyle \,{\cal L  }_t\,\,H\, \,\,= \frac{1}{\sqrt{D}} \, \frac{e^{- \frac{|x|}{\sqrt{D}}\,\sigma }}{2 \,\sigma \,  }
\end{equation}

\noindent where 
\begin{eqnarray}
  \displaystyle   \sigma^2 \ \,=\, s\, +\, a \, + \, \frac{\delta}{s+\delta d}\, + \, \frac{\varepsilon}{s+\,\beta \varepsilon}.\,\,
\end{eqnarray}

Moreover, function H(x,t)  satisfies some  properties typical of the  fundamental solution of heat equation, such as:

 a) \hspace{1mm} $ H( x,t) \, \, \in  C ^ {\infty}, \,\,\,\,$ \,$\,\,\, t>0, \,\,\,\, x \,\,\, \in \Re, $ 

b) \hspace{1mm} for fixed $\, t\,>\,0,\,\,\, H \,$  and its derivatives are vanishing exponentially fast as $\, |x| \, $  tends to infinity.

c) \hspace{1mm} In addition, it results   $ \displaystyle\lim _{t\, \rightarrow 0}\,\,H(x,t)\,=\,0, $ for any fixed $\, \eta \,>\, 0,\, $ uniformly for all $\, |x| \,\geq \, \eta.\, $ \hbox{}\hfill\rule{1.85mm}{2.82mm}

 \vspace{3mm} To obtain results of existence and uniqueness for the problem (\ref{19}), the theorem of fixed point can be applied  and therefore, also according to \cite{c}, for initial term and source function  we shall admit:

\vspace{3mm}
   
\noindent   \textbf{Assumpion A} Initial data $ u_0  $ is continuously differentiable and bounded together with its first derivative. 
 The source term    
     $ F(x,t,u)$   is defined and continuous on the following set: 
     
     \begin{equation}  \label {z}
     Z = \{ (x,t,u ) : (x,t) \in \Omega_T , -\infty<u<\infty \}.
     \end{equation}

 Besides,  for each $ K>0  $  and $ |u| < K, $ $\, F(\,x,\,t,\,u\,)\,$ is uniformly Lipschitz 
 continuous in $(x,t) $  for each compact set of   $ \Omega_T $ and it is bounded for bounded $ u.   $   
 
 Then,  for all ($ u_1,u_2 ),$  there exists a positive  constant $ W_F $  such that: 
 
 \begin{equation}
      | F(\,x,\,t,\,u_1\,) -F(\,x,\,t,\,u_2) | \leq W_F \,\,|u_1-u_2|.
     \end{equation}\hbox{}\hfill\rule{1.85mm}{2.82mm}

\vspace{1mm}  
As a consequence, when  the fundamental solution  $ H(x,t)  $ and source function $ F(x,t,u) $  satisfy theorem \ref{th} and Assumption A, respectively, indicating by $ u(x,t) $  
a solution of problem $\,{\cal P}, $  then 
$ u $  assumes this form:

\begin{eqnarray}  \label{A14}
 & \displaystyle \nonumber u (x,t)  =\int_\Re   \,H ( x-\xi, t)\,\, u_0 (\xi)\,\,d\xi \,
 \\
 \\
 &\displaystyle \nonumber\,+\,\int ^t_0     d\tau \int_\Re   H ( x-\xi, t-\tau)\,\, F\,[\,\xi,\tau, u(\xi,\tau\,)\,]\,\, d\xi.
\end{eqnarray} 

\vspace{3mm} \noindent On the other hand,   if $ u(x,t)   $ is a continuous  and bounded solution of (\ref{A14}), it is possible to prove that $ u $  satisfies  (\ref{19}).

\vspace{3mm} \noindent Consequently, it is possible to conclude that 

\begin{theorem}
Initial  value problem  (\ref{19})  admits a unique  solution  only if   $(\ref {A14})$  admits a unique continuous and bounded solution. 
\hbox{}\hfill\rule{1.85mm}{2.82mm}
\end{theorem}

Besides, by means of  fixed point theorem,(and  extensive proofs can be found, f.i., in \cite{c,32,df13,ddf,dr8}), it is possible to prove the following theorem: 

\begin{theorem}
When Assumption A is satisfied, then the initial value problem (\ref{19}) admits a unique regular solution $ u(x,t)  $ in $ \Omega_T. $
\hbox{}\hfill\rule{1.85mm}{2.82mm} 
\end{theorem}

In this case, taking into account  the  source term $ F(x,t) $ defined in (\ref{18}), solution  (\ref{A14}) assumes the following form:

\begin{eqnarray}  \label{abc}
 \nonumber & \displaystyle u(x,t) = \int ^t_0     d\tau \int_\Re   H ( x-\xi, t-\tau) \varphi\,[\xi,\tau, u(\xi,\tau)] d\xi 
 \\\nonumber
\\ 
&\nonumber\displaystyle   + \bigg (\frac{h}{d}- \frac{c}{\beta}\bigg) \int ^t_0     d\tau \int_\Re   H ( x-\xi, t-\tau)  d\xi +  \frac{c}{\beta}  \int ^t_0    e^{-\beta \varepsilon  \tau }  d\tau \int_\Re    H ( x-\xi, t-\tau)d\xi  
\\  
\\
&\nonumber\displaystyle -  \int ^t_0  e^{-\beta \varepsilon  \tau }    d\tau \int_\Re    H ( x-\xi, t-\tau) w_0(\xi) d \xi  - \frac{h}{d} \int ^t_0    e^{-\delta d  \tau }   d\tau \int_\Re   H ( x-\xi, t-\tau)d\xi
\\\nonumber
\\\nonumber
&\nonumber\displaystyle + \int ^t_0  e^{-\delta d \tau }   d\tau   \int_\Re      H ( x-\xi, t-\tau)  y_0(\xi) d\xi + \int_\Re   H ( x-\xi, t) u_0 (\xi)\,\,d\xi 
\end{eqnarray}

\vspace{3mm}\noindent and   this formula, together with  relations (\ref{17}), allows us  to determine  also  $\, v(x,t) \,$  and $\, y(x,t) \,$  in terms of the data.

\section{Some properties related to H(x,t)}

In order to  obtain a priori estimates and asymptotic effects,
some properties related to the fundamental solution $ H $ need to be  evaluated.

More precisely, formula (\ref{abc}) shows the need
 to evaluate  the convolution of the fundamental solution  $ H $ with respect to  time and space.
 
Consequently, this section will include a first part where  two theorems involving some properties related to $ H(x,t) $ are showed,  and a subsection where some  premises allowing  to prove properties related to convolution integrals, will be stated.

\vspace{3mm}Let us  start  indicating  by

\begin{equation} \label{a218}
A(t) =\frac{e^{-\beta \varepsilon \, t}- e^{- a t}}{a- \beta \varepsilon }; \quad B(t) =\frac{e^{- \delta d  t}-e^{-a t} }{a-\delta d  };  \quad C(t) =\frac{ e^{- \delta d  t}-e^{-\beta \varepsilon \, t}}{\beta \varepsilon -\delta d  } 
\end{equation}

\noindent three positive functions, then  the following theorem holds:  
\begin{theorem}
The solution function $ H  $ defined in (\ref{A16}) satisfies the following estimate:

\begin{equation} \label{a1}
| H| \leq  \frac{e^{- \frac{x^2}{4\,D\, t}\,}}{2 \sqrt{\pi  D t } }\,\,\,\bigg[ e^{-\,a\,t}+ t  \varepsilon A(t) + \delta \, t  \bigg(1+ \frac{\varepsilon t }{|a-\beta \varepsilon|} \bigg)B(t)\,+ \frac{\varepsilon t}{|a-\varepsilon \beta|}\, C(t)\bigg]
\end{equation}
\end{theorem}


Since

  \begin{equation}
  |J_1(\,2 \,\sqrt{\,\varepsilon \,y\,(t-y)\,}\,\,) |  \leq \sqrt{\,\varepsilon \,y\,(t-y)\,}\,\,\ \quad  (y\leq t)
  \end{equation}

\noindent from (\ref{22}) it results:

\begin{eqnarray} 
&\nonumber \displaystyle | H_1(x,t)|\,\leq\, \, \, \frac{e^{- \frac{x^2}{4\,D\, t}\,}}{2 \sqrt{\pi  D t } }\,\,\bigg[\, e^{-\,a\,t}\,+ \varepsilon \,t\,\,
 \int^t_0  e^{- a\,y} \,\, e^{-\beta \varepsilon \,(\, t \,-\,y\,)} \, dy\bigg]
 \end{eqnarray}
 
 \noindent and hence:

\begin{equation}  \label {222}
|H_1(x,t)|\,\leq\, \, \, \frac{e^{- \frac{x^2}{4\,D\, t}\,}}{2 \sqrt{\pi  D t } }\,\,\bigg[\, e^{-\,a\,t}\,+ \varepsilon \,t\, \, \frac{e^{-\beta \varepsilon \,\, t \,\,}  - e^{- at}}{a-\varepsilon \beta}\bigg].
\end{equation}

 \vspace{3mm} Moreover, from  (\ref{H2}) and by means of (\ref{222}), it results:

\begin{eqnarray}  
 & \nonumber \displaystyle | H_2 |\leq \int _0 ^t  \frac{e^{- \frac{x^2}{4\,D\, y}\,}}{2 \sqrt{\pi  D y } }\,\,\,\bigg[ e^{-\,a\,y} +
   \varepsilon \,y\,\,\,   \,\,\,\frac{e^{-\beta \varepsilon \,\, y \,\,}  - e^{- ay}}{a- \varepsilon \beta} \bigg] \,e^{ -\delta d (t-y)} \,\,\delta y   \,\, dy.  
 \end{eqnarray}

\vspace{3mm} Consequently one obtains:

 \begin{equation} \label{aaa}
   |H_2| \leq      \frac{  \delta t  \,\,e^{- \frac{x^2}{4\,D\, t}\,}}{2 \sqrt{\pi  D t } } \bigg[\frac{ e^{-\delta d t}- e^{-a\,t}}{ a-\delta d }    \bigg(1 + \frac{\varepsilon   \,t}{|a- \varepsilon \beta|}  \bigg) + \frac{\varepsilon   \,t}{|a- \varepsilon \beta|} 
   \frac{ e^{-\,\delta d t}-e^{-\beta \varepsilon}}{\beta \varepsilon-\delta d}\bigg]   
    \end{equation}
    
\vspace{3mm}    Hence, according to (\ref {A16}), for (\ref{222})  and (\ref{aaa}),  theorem  holds. \hbox{}\hfill\rule{1.85mm}{2.82mm}


Now, let us introduce  as $ I_0  $   the modified Bessel function of the first kind  and order $ 0, $ and let

\begin{equation}\label{e46}
l= \min (a, \beta \varepsilon),\qquad q= \min \{ a, \beta \varepsilon , \delta d \},
\end{equation}

\begin{equation} \label{e47}
\lambda (t)\equiv  1+  \pi t (\sqrt{\varepsilon } +\sqrt{\delta } + \pi t \sqrt{\delta \varepsilon  }) .\end{equation} 
 
\vspace{2mm} The following theorem holds:

\begin{theorem}

The fundamental solution $ H(x,t)  $  defined in (\ref{A16}) satisfies the following estimates:

\begin{eqnarray}    \label{pr2} 
&\nonumber \displaystyle \int_ \Re   | H(x-\xi,t) | d\xi \leq e^{-at} +\sqrt{\varepsilon } \pi \,t \,\,e^{- \frac{\beta \varepsilon +a}{2}t}  \,\, I_0 \bigg(\frac{\beta \varepsilon -a }{2} t \bigg) 
\\   
\\
&\nonumber +\,\sqrt{\delta} \,     \pi \,t \, \bigg[\,e^{- \frac{\delta d  +a}{2}t}  \,\, I_0 \bigg(\frac{\delta d  -a }{2} t \bigg)  +  \sqrt{\varepsilon } \pi \,t  e^{- \frac{\delta d +l}{2}t}  \,\, I_0 \bigg(\frac{\delta d  -l }{2} t \bigg) \bigg];
\end{eqnarray}

\begin{eqnarray}  \label{439}
&\displaystyle \int_ \Re    |H (x-\xi,t)|  d\xi \leq    \lambda (t)\,e^{-q t}
\end{eqnarray}

Besides, indicating by

\begin{eqnarray} \label{2s}
& \displaystyle S= 1/a+ \sqrt{\varepsilon } \,\pi  \,\,\frac{a+\beta \varepsilon }{2 (a \beta \varepsilon ) ^{3/2}}+ \sqrt{\delta } \,\pi \bigg[ \,\,\frac{\delta d +a}{ (a \delta d  )   ^{3/2}} + 3\pi \sqrt{\varepsilon }\,\,\frac{\delta^2 d^2 +l^2}{4 (l \delta d  )  ^{5/2}}\bigg],  
\end{eqnarray}

\noindent one has:

\begin{eqnarray}  \label{a441}
&\displaystyle \int _0^t\, d\tau\int_ \Re   | H(x-\xi,t-\tau) | d\xi \leq   S.  
\end{eqnarray}

\end{theorem}

Considering that 
 
\begin{equation} \label{bbb}
\displaystyle H = H_1  - H_2,
\end{equation}
\noindent we will  firstly  focus on the integral involving $ H_1, $ and  then on that involving $ H_2. $

 \vspace{3mm}Since it results:  

\begin{equation}
\int_ \Re  e^{-\frac{x^2}{4Dt}} d x = 2 \sqrt{\pi Dt} ;\qquad  \qquad  |J_1(z)| \leq 1,
\end{equation} 

\vspace {3mm} \noindent  from (\ref{22}) one obtains:

\begin{eqnarray}  \label{u}
& \displaystyle  \int_ \Re  |H_1(x,t)| dx  \leq   e^{-\,a\,t}\,
 \displaystyle + \sqrt{\varepsilon} \int^t_0   e^{-\beta \varepsilon (t-y)  }\,\,e^{-ay} \frac{\sqrt{y}}{\sqrt{t-y}} \,\,  dy 
\end{eqnarray}

\vspace{3mm}\noindent with
\begin{eqnarray}
&\nonumber\displaystyle \int ^t_0   e^{-\beta \varepsilon (t-y)  } e^{-ay} \sqrt{\frac {y}{t-y} }  dy  =
  -\int^t_0  e^{-\beta \varepsilon (t-y)  }e^{-ay}  (t/2-y)\frac{dy}{\sqrt{y(t-y)} }+  
  \\\nonumber\displaystyle
  \\
  \\\nonumber\displaystyle
 & \nonumber\displaystyle + \int^t_0 e^{-\beta \varepsilon (t-y)  }e^{-ay} \frac{t/2 \,dy }{\sqrt{y(t-y}) }.   
\end{eqnarray}

\vspace{3mm}\noindent Now, taking into account that

\begin{equation}
      \displaystyle\int _0^{2b} e^{-sy }  \,  (b-y) \frac{1}{\sqrt{2by-y^2}} dy = \pi b e^{-sb} I_1 ( sb)    
\end{equation}

\vspace{3mm}\noindent and

\begin{equation}
     \displaystyle\int _0^{2b} e^{-sy }  \,   \frac{1}{\sqrt{2by-y^2}} dy = \pi  e^{-sb} I_0 ( sb)     
\end{equation}

\vspace{3mm}\noindent for $ b=t/2$  and $ s= a-\beta \varepsilon, $  one has:

\begin{eqnarray}
&\nonumber\displaystyle \int ^t_0   e^{-\beta \varepsilon (t-y)  } e^{-ay} \sqrt{\frac {y}{t-y} }  dy      =
  \frac{ \pi \,t}{2} \, \bigg[  e^{- \frac{a-\beta \varepsilon }{2}t}      \bigg(I_0( \frac{a-\beta \varepsilon}{2}t) -  I_1(\frac{a-\beta \varepsilon}{2}t )\bigg) \bigg].
\end{eqnarray}

\vspace{3mm} Consequently, as for

\vspace{3mm}\begin{center}
$ I_1(-z)=-I_1 (z) \qquad  I_0(z)=I_0(-z) \qquad  I_1(|z|)\leq I_0 (|z|),  $ 
\end{center}
\noindent   it results:

\begin{equation} \label{443}
\displaystyle \int_ \Re    |H_1(x-\xi,t)|  d\xi \leq  e^{-at} +\sqrt{\varepsilon } \pi \,t \,\,  e^{- \frac{ a+\beta \varepsilon }{2}t}  \, \,\, I_0 \bigg(\frac{\beta \varepsilon -a }{2} t \bigg).
\end{equation}

\vspace{3mm} Now,  being    $  I_0(|z| )< e^{|z|}$, from (\ref{443}) one deduces that

\begin{equation}   \label{445}
 \displaystyle \int_ \Re    |H_1 (x-\xi,t)|  d\xi \leq  e^{-at} +\sqrt{\varepsilon } \pi \,t \,\, \,\, e^{- l t}
\end{equation}

 \vspace{3mm}\noindent where  $ l$ is defined in $ (\ref{e46})_1. $

\vspace{2mm}As for function $ H_2, $ taking into account that  $ |J_1|\leq  1,$  from (\ref{H2}) and  by means of (\ref{445}), it results:

\begin{eqnarray} 
&\nonumber \displaystyle \int_ \Re| H_2 (x-\xi,t)| \, d\xi  \leq \,\sqrt{\delta} \int _0 ^t  \big(e^{-ay} +\sqrt{\varepsilon } \pi \,y \,\ e^{- l y}\big) \,\,e^{ -\delta d (t-y)} \sqrt{\frac{ y}{t-y}}   \, dy.
\end{eqnarray}

\vspace{3mm} Hence, returning to the previous reasoning, one obtains:

\begin{equation}  \label{446}
\displaystyle \int_ \Re    |H_2| \leq  \,   \sqrt{\delta}\,  \pi \,t \, \bigg[ \,e^{- \frac{\delta d  +a}{2}t}   I_0 \bigg(\frac{\delta d  -a }{2} t \bigg)  +  \sqrt{\varepsilon\, }\, \pi \, t \, e^{- \frac{\delta d +l}{2}t}   I_0 \bigg(\frac{\delta d  -l }{2} t \bigg) \bigg]
\end{equation}

\noindent from which, along with (\ref{443}), (\ref{pr2})   is proved.

\vspace{3mm} Moreover, from  (\ref{446}), an inequality analogous to (\ref{445}) can be obtained. In this way,  according to (\ref{bbb}),  (\ref{439}) follows, too.

\vspace{2mm}Lastly, since it results

\begin{equation}
\int_0^\infty e^{-pt }\, t\, I_0 (bt) \, dt = p\, (\sqrt{p^2-b^2})^{-3} \qquad  Re\,\, p > |Re \,\,b|
\end{equation}

\begin{equation}
\int_0^\infty e^{-pt }\, t^2\, I_0 (bt) \, dt =  (\sqrt{p^2-b^2})^{-3/2} \bigg(\frac{3 p^2}{p^2-b^2} -1 \bigg) \qquad  Re\,\, p > |Re \,\,b|,
\end{equation}

\vspace{3mm} \noindent from  (\ref{443}) and (\ref{446}), property  (\ref{a441}) can be proved. 
\hbox{}\hfill\rule{1.85mm}{2.82mm}

\subsection{Premises on convolution integrals referring to the solution }

In order to determine the estimates related to the solution, it is necessary to highlight every convolution integrals that characterize the solution itself. Therefore, in this subsection convolutions  $ K_\delta $  and  $  H_\delta $   will be introduced and, by means of them, solution $ u(x,t)$ will be expressed.(Formula (\ref{318})).

\vspace{2mm}Hence,  
 let us  consider

\begin{eqnarray} \label{c38}
 K_ \delta (x,t) \equiv \int^t_0 \, e^{- \delta d \,(t- y)}\,  H_1 (x,y)  \, J_0 \, (\, 2\, \sqrt{\delta\, y ( t-y)  }\,) \,\,dy 
\end{eqnarray}

\noindent  and let

\begin{equation}
\,\, g_1(x,t)\,\ast g_2(x,t)  = \int _0^t g_1(x,t-\tau) g_2(x,\tau) \,d\tau
\end{equation} 
be the convolution with respect to $ t. $

\vspace{2mm} In \cite{2020} it has been proved that:

\begin{equation}                              \label{317}
  e^{-\,\delta d \,t} \ast\,H =   K_\delta
\end{equation}

\noindent and

\begin{eqnarray} \label{A41}
e^{-\,\varepsilon\,\beta \,t} *\,H=   K_\delta+ (\delta d-\varepsilon \beta ) e^{-\beta \varepsilon  \,t} * K_\delta.
\end{eqnarray}

\vspace{3mm} Now, denoting by

\begin{equation} \label{ss}
 H_\delta =   \int^t_0 e^{- \varepsilon \beta \,(t-\tau)}  d\tau \int^\tau_0   H_1 \,(x,y)\, e^{-\delta d (\tau-y)}     J_0 (\, 2\, \sqrt{\delta y ( \tau-y)  }) dy  
\end{equation}

\vspace{3mm}\noindent it results:

\begin{equation}                              \label{aa317}
  H_\delta \equiv e^{-\, \beta \varepsilon  \,t} \ast\,K_\delta, 
\end{equation}

\vspace{2mm} 
\noindent and as a consequence,  from (\ref{A41}), one one: 

\begin{eqnarray} \label{xxx}
e^{-\,\varepsilon\,\beta \,t} *\,H=   K_\delta+ (\delta d-\varepsilon \beta )\, H_\delta.
\end{eqnarray}

\vspace{2mm}  Moreover, let us denote  by 
  
\begin{equation}
\,\, g_1(x,t)\,\diamondsuit  \, g_2(x,t)  = \int _{\Re} f_1(\xi,t) g_2(x-\xi,t) \,\,  d\xi
\end{equation}   
  
\vspace{3mm} \noindent the convolution with respect to  the space, and   

\begin{equation}
H  \otimes F\,=\, \int_0^t\,d\tau\, \int_\Re \,H(x-\xi,t-\tau) \, \,F \,[\,\xi, \tau ,u(\xi,\tau)\,]\,d\xi.
\end{equation}

\vspace{3mm}  Since (\ref{317}) and (\ref{xxx}),  it results:

\vspace{3mm}

\begin{equation}
\label{A51}
\left \{
   \begin{array}{lll}
\displaystyle H  \otimes \, e^{-\delta d t }=\int_\Re   K_\delta (\xi,t)\,\,d\xi,
\\
\\ 
 \displaystyle H  \otimes \, e^{- \beta \varepsilon\, t } =\int_\Re \big[  K_\delta+ (\delta d-\varepsilon \beta ) H_\delta\,\big]  d\xi
  \end{array}
  \right.  
  \end{equation}
  
\vspace{3mm} \noindent and

\begin{equation}
\label{315}
\left \{
   \begin{array}{lll}
H  \otimes  ( y_0(x)\, e^{-\delta d t })= y_0  \,\diamondsuit \, K_\delta
  \\
  \\
H  \otimes  ( w_0(x) \,e^{-\beta \varepsilon t })= w_0 \diamondsuit  [  K_\delta+( \delta d-\varepsilon \beta )  H_\delta\,].  
 \end{array}
  \right.  
\end{equation}

\vspace{3mm} Consequently, given  (\ref{abc}) , we get:

\begin{eqnarray}  \label{318}
 \nonumber & \displaystyle u(x,t) \,=\, H  \,\diamondsuit \, u_0 (x) \,  + K_\delta \,\diamondsuit \, (y_0(x) -\, w_0(x))    + \, \,  H \,  \otimes  \varphi (u)  \, \,
\\ \nonumber
\\ 
&\nonumber\displaystyle  +(\varepsilon \beta - \delta d)\,H_\delta \,\diamondsuit\, w_0(x) \, +  \frac{c}{\beta}\, H_\delta   \,\diamondsuit \,  \big ( \delta d - \varepsilon \beta )\,   
\\
\\
\nonumber &\displaystyle  + H  \,\otimes \,\bigg (\frac{h}{d}\,- \frac{c}{\beta}\,\bigg)  +  K_\delta \,\diamondsuit \,    \bigg (\frac{c}{\beta}-\frac{h}{d}\,\bigg)
\end{eqnarray}

\vspace{3mm}\noindent and   this formula explicitly shows all the convolutions involved in the solution $ u(x,t). $

\section{ On convolutions involving functions $ K_\delta$ and $ H_\delta  $}

Formula  (\ref{318}) shows that an analysis  of the  solution  directly implies estimates   on  both    $ H(x,t)  $ and   on  functions $ K_\delta,  $  $ H_\delta,  $ defined in  (\ref{c38}) and (\ref{ss}).  

For this, let us consider $ A(t), B(t), C(t), \lambda(t)  $   defined in (\ref{a218}) and (\ref{e47}), respectively. Moreover, let

 \begin{equation}  \label{hEE}
E(t) = \frac{ e^{-qt}-e^{- \delta d  t}}{\delta d -q  }   \qquad L(t)=\frac{ e^{-qt}-e^{- \beta \varepsilon   t}}{\beta \varepsilon -q  } 
\end{equation}

\vspace{3mm} \noindent with $ q $  defined by $(\ref{e46})_2$\,.

\vspace{3mm} In addition,  

\begin{equation} \label{E}
 M= \frac{1}{|\delta d-q| \delta q d} \bigg[ q + \delta d + \pi (\sqrt{\varepsilon} +\sqrt{\delta }) \big(\frac{q^2+\delta^2 d^2}{\delta d q} \bigg)+ 2\pi^2 \sqrt{\delta \varepsilon} \bigg( \frac{q^3 + \delta^3 d^3}{(q \delta d)^2}\bigg) \bigg]
\end{equation}

\begin{equation} \label{EE}
 N= \frac{1}{|\beta \varepsilon-q| q\beta \varepsilon} \bigg[ q + \beta \varepsilon + \pi (\sqrt{\varepsilon} +\sqrt{\delta }) \big(\frac{q^2+\beta^2 \varepsilon^2 }{\beta \varepsilon  q} \bigg)+ 2\pi^2 \sqrt{\delta \varepsilon} \bigg( \frac{q^3 + \beta ^3 \varepsilon^3}{(q \beta d)^2}\bigg) \bigg]
\end{equation}

\begin{equation} \label{g}
g(t)=    \frac{\lambda (t) }{|\beta \varepsilon -\delta d|} \,  \big[ E(t) +    L(t)  \big]   
\end{equation}

\begin{equation} \label{h}
 h(t) = \frac{  \, \lambda (t) }{( \varepsilon \beta -\delta d )^2}   \big[L(t) + (1+t(\delta d - \varepsilon \beta  \big] E(t).
 \end{equation}

\vspace{3mm} \noindent The following theorems hold:

 \begin{theorem}  \label{uff}
 Function  $ K_\delta (x,t) $  defined in  (\ref{c38}) satisfies the following estimates:
 
 \vspace{3mm}
\begin{equation} \label{b38}
\int_ \Re \big| K_ \delta (x,t)\big | \leq \, \lambda (t) \,E(t);  
\end{equation}

\vspace{3mm} 
\begin{equation}  \label{bb38}
\qquad \int _0^t d\tau \int_ \Re \big| K_ \delta (x,\tau) \big| dx \leq  M.
\end{equation}

\vspace{3mm}
\begin{equation} \label{dccc}
\int _0^t   e^{-\delta d \tau}  d\tau\int_ \Re \big| K_ \delta (x,t-\tau) \big| dx \leq      t \, \lambda(t)\,  E(t)
\end{equation}
 \end{theorem}

\vspace{3mm}
By means of  (\ref{317}) and  property (\ref{439}) on $ \int_\Re |H(\xi,t)| d\xi\, $, inequality (\ref{b38}) follows. 

\vspace{3mm}\noindent By this estimate, according to (\ref{e47}), and taking into account that 

\begin{equation}
\int_ 0^t y \, e^{-\alpha y} \leq 1/\alpha^2; \qquad \int_ 0^t y^2 \, e^{-\alpha y} \leq 2/\alpha^3  \qquad (t>0, \quad \alpha >0),
\end{equation}

\vspace{3mm}\noindent  (\ref{bb38}) holds, too.

\vspace{3mm}\noindent  Moreover, because of  (\ref{317}), it results

\begin{equation}
e^{-\delta d t }\ast \,K_\delta = e^{-\delta d  t}\ast \,H \ast e^{-\delta d t} = (t \,\,e^{- \delta d t})  \ast  H  
\end{equation}

\noindent and inequality  (\ref{dccc}) follows.
\hbox{}\hfill\rule{1.85mm}{2.82mm}

 \begin{theorem}
Referring to (\ref{ss}), function $ H_\delta(x,t)  $  satisfies the  inequalities below:
 \end{theorem}

\begin{equation} \label{bbb38}
 \int_ \Re | H_ \delta (x,t) |\,dx\,\,\leq  g(t)
\end{equation}

\begin{equation}\label{ddccc}
\int _0^t    d\tau\int_ \Re \big| H_ \delta (x,t-\tau)\big| dx \leq \,\frac{M+N}{|\beta \varepsilon-\delta d|}
\end{equation}

\begin{equation} \label{bccc}
\int _0^t   e^{-\delta  d \tau}  d\tau\int_ \Re \big| H_ \delta (x,t-\tau)\big| dx \leq   h(t).
\end{equation}

\begin{equation} \label{ccc}
\int _0^t   e^{-\beta \varepsilon \tau}  d\tau\int_ \Re \big| H_ \delta (x,t-\tau) \big| dx \, \leq \, \frac{  t \,\, \lambda (t)}{|\delta d-q| } \,\,\big[ C(t) + L(t)\big]  
\end{equation}

\vspace{3mm} 

According to (\ref{317}) and (\ref{aa317}), one has:

\begin{equation}
\int_ \Re | H_ \delta (x,t) |\,dx\, = \int _0^t   e^{-\beta \varepsilon \tau} \, d\tau\int_ \Re  \big |K_ \delta (x,t-\tau)\big| dx
\end{equation}

\noindent with 

\begin{equation}
e^{-\beta \varepsilon t}\ast \,K_\delta = e^{-\beta \varepsilon t}\ast \,H \, \ast e ^{-\delta d  t}  = C(t) \ast  H(x,t)  
\end{equation}

\vspace{3mm} \noindent  where $ C(t)  $  is  defined in $(\ref{a218})_3.$ Hence, since (\ref{439}), inequality (\ref{bbb38}) holds.
 
\vspace{3mm} 
\noindent Consequently, also (\ref{ddccc}) follows.

\vspace{3mm} \noindent  Estimate (\ref{bccc}) is proved  by means of  

\begin{equation}
e^{-\delta d t }\ast \,H_\delta = e^{-\delta d  t} \ast \,K_\delta  \ast e^{-\beta \varepsilon  t} = (t \,\,e^{- \delta d t})  \ast   e^{-\beta \varepsilon  t} \ast  H.   
\end{equation}

Finally,  taking into account that

\begin{equation}
e^{-\beta \varepsilon t}\ast \,H_\delta = e^{-\beta \varepsilon t}\ast \,K_\delta \ast e^{-\beta \varepsilon t} = (t \,\,e^{-\beta \varepsilon t})  \ast  K_\delta,   
\end{equation}

\noindent  from  (\ref{b38}),  (\ref{ccc}) is proved,too.
\hbox{}\hfill\rule{1.85mm}{2.82mm}

\section{Analysis of solution}

In order to analyse  functions $u(x,t), w(x,t),$ and $y(x,t),$ it  appears necessary to make explicit the integrals of convolutions involving  functions
$H_\delta$ and $ K_\delta $ whose estimates have been established in the previous section.

Therefore, since   (\ref{318}), by means of convolution properties, we get:

\begin{eqnarray}  \label{abcc}
 \nonumber & \displaystyle u(x,t) = \int ^t_0     d\tau \int_\Re   H ( x-\xi, t-\tau) \varphi\,[\xi,\tau, u(\xi,\tau)] d\xi 
 \\\nonumber
\\ 
&\nonumber\displaystyle   + \bigg (\frac{h}{d}- \frac{c}{\beta}\bigg) \bigg[ \int ^t_0     d\tau \int_\Re   H ( x-\xi, t-\tau)d\xi  -  \int_\Re   K_\delta ( x-\xi, t)d\xi\bigg]
\\ 
\\
&\nonumber\displaystyle   +   \int_\Re      K_\delta ( x-\xi, t)  \big[ y_0(\xi) -  w_0(\xi)\big] d\xi - (\delta d-\varepsilon \beta )  \int_\Re    H_\delta( x-\xi, t)\, w_0(\xi) d \xi
\\\nonumber
\\\nonumber
&\nonumber\displaystyle   +\frac{c}{\beta}  (\delta d-\varepsilon \beta ) \int_\Re    H_\delta  ( x-\xi, t)d\xi  + \int_\Re   H ( x-\xi, t)\, u_0 (\xi)\,\,d\xi. 
\end{eqnarray}

\vspace{5mm} Moreover, as for functions  $ w(x,t) $  and $ y(x,t)$ defined in (\ref{17}), according to (\ref{317}), (\ref{aa317}) and (\ref{xxx}), since (\ref{abcc}), the following integrals must be considered:

\begin{eqnarray}  \label{yabcc}
 \nonumber & \displaystyle \int_0 ^t e^{-\beta \varepsilon (t-\tau) }   u(x,\tau)d \tau  =  \int_\Re   K_\delta  ( x-\xi, t) u_0 (\xi)\,\,d\xi\
 \\\nonumber
\\ 
&\nonumber\displaystyle +\, \int ^t_0     d\tau \int_\Re   K_\delta ( x-\xi, t-\tau) \bigg[\varphi\,[\xi,\tau, u(\xi,\tau)] + \frac{h}{d}- \frac{c}{\beta} \bigg]\,d\xi  
\\
\\
\nonumber & \displaystyle  + (\delta d-\varepsilon \beta )  \int ^t_0     d\tau \int_\Re   H_\delta ( x-\xi, t-\tau) \big[\varphi\,[\xi,\tau, u(\xi,\tau)] +\frac{h}{d}- \frac{c}{\beta}\big] d\xi 
\\\nonumber
\\\nonumber
&\nonumber\displaystyle  + (\delta d-\varepsilon \beta ) \int_0 ^t e^{-\beta \varepsilon  (t-\tau) } d\tau \int_\Re    H_\delta( x-\xi, t-\tau) \big[\frac{c}{\beta} - w_0(\xi) \big ] \,\,d \xi
\\\nonumber
\\ 
&\nonumber\displaystyle   +   \int_\Re      H_\delta ( x-\xi, t-\tau) \big[  y_0(\xi) -  w_0(\xi) -\frac{h}{d}+ \frac{c}{\beta} + (\delta d-\varepsilon \beta ) u_0(\xi)\big] d\xi  
\end{eqnarray}

\vspace{5mm} \noindent  and

\begin{eqnarray}  \label{yabccc}
 \nonumber & \displaystyle \int_0 ^t e^{-\delta d (t-\tau) }   u(x,\tau)d \tau  =  \int_\Re   K_\delta  ( x-\xi, t) u_0 (\xi)\,\,d\xi
 \\\nonumber
\\\nonumber
&\nonumber\displaystyle  +\int ^t_0     d\tau \int_\Re   K_\delta ( x-\xi, t-\tau) \bigg[\varphi\,[\xi,\tau, u(\xi,\tau)] + \frac{h}{d}- \frac{c}{\beta} \bigg]\,d\xi  
\\
\\\nonumber
&\nonumber\displaystyle  + \int_0^t e^{-\delta d t } d\tau\int_\Re      K_\delta ( x-\xi, t-\tau) \big[  y_0(\xi) - w_0(\xi)- \frac{h}{d}+ \frac{c}{\beta}\, \big] d\xi
\\\nonumber
\\\nonumber
&\nonumber\displaystyle  + (\delta d-\varepsilon \beta ) \int_0 ^t e^{-\delta d  (t-\tau) } d\tau \int_\Re    H_\delta( x-\xi, t-\tau) \big[\frac{c}{\beta} -  w_0(\xi) \big ] \,\,d \xi.
\end{eqnarray}

\section{Estimates of solution}

As for the analysis  of solutions of the non linear reaction diffusion model, there exists a large bibliography. In particular in \cite {smoller,R21} the existence of bounded solutions  is proved.

\noindent Therefore, in the class of bounded solutions, let us
assume initial data and  function  $ \varphi(x,t,u)  $ satisfy Assumption A, and let

\vspace{5mm}

\[ ||\,u_0\,|| \,= \displaystyle \sup_ { \Re}\, | \,u_0 \,(\,x\,) \,|, \quad ||\,w_0\,|| \,= \displaystyle \sup _{\Re }\, | \,w_0 \,(\,x\,) \,|,  \,\quad ||\,y_0\,|| \,= \displaystyle \sup _{\Re }\, | \,y_0 \,(\,x\,) \,|, \]

 \[ ||u|| =  \displaystyle \sup _{ \Omega _T\,} | \,u(x,t)\, \qquad
||\varphi|| \,= \displaystyle \sup _{ Z\,} | \,\varphi \,(\,x,\,t,\,u) \,| \]

\vspace{2mm} \noindent  with  $ \varphi$ defined  in  $(\ref{14})$  and $ Z  $  defined in  (\ref{z}).

\vspace{1mm}In order to give a priori estimates of the solution of FHR system, the following theorem is proved: 
 
 \begin{theorem}
 If function  $\varphi(x,t,u) $ and initial data $ u_o(x),\, w_o(x),\, y_o(x) $ are compatible with Assumption A, then the problem (\ref{11})-(\ref{ic}) satisfies the following estimates: 
 
 \begin{equation}
\label{68}
     \begin{array}{lll}
    \displaystyle |u (x,t) |& \leq  ||u_0 (x)|| \,\,\displaystyle  \lambda (t)  \,\displaystyle  e^{-q t}+ \bigg (|| \varphi|| + \displaystyle \bigg|\frac{h}{d}- \frac{c}{\beta}\,\bigg| \bigg ) \,S  
\\ 
\\
&\displaystyle + \bigg( \displaystyle ||y_0||+||w_0)||+\bigg|\frac{h}{d}- \frac{c}{\beta}\,\bigg|\bigg)   \,\,\lambda (t)\,\, E(t)
\\ 
\\
&\displaystyle +  \bigg( ||w_0||+\frac{c}{\beta}    \bigg)  \displaystyle (|\delta d-\varepsilon \beta |) \,\, g(t);
\end{array}
\end{equation}

\begin{equation}   \label{fff}
 \begin{array}{lll}
 \displaystyle  |w (x,t) |& \leq  \displaystyle ||w_0||  \,\,e^{-\beta \varepsilon t} +\displaystyle  \frac{c}{\beta}  + \varepsilon \,\,| |u_0|| \,\,\lambda (t)\,\, E(t)
\\ \\
&\displaystyle + \varepsilon \bigg( ||\varphi||+ \bigg|\frac{h}{d} -  \frac{c}{\beta}\bigg|    \bigg)  (2M+N) +   
 \\
 \\
 &\displaystyle +\varepsilon  \frac{|\delta d - \varepsilon \beta| }{|\delta d-q| } \,\bigg[\,\frac{c}{\beta} + ||w_0(x)|| \,\bigg  ] t \,\lambda (t)\, \big[C(t)\, + \,L(t) \big]
\\
\\
&\displaystyle +\varepsilon \bigg[  ||y_0||+  ||w_0|| +\bigg|\frac{c}{\beta}-\frac{h}{d}\bigg|+  |\delta d-\varepsilon \beta | \,\, ||u_0||\bigg] g(t);
\end{array} 
\end{equation}\\

\begin{equation}  \label{ff}
 \begin{array}{lll}\displaystyle  |y (x,t) | &\leq \displaystyle ||y_0||  e^{\,-\,\delta d t} \,+\, \frac{h}{d}\,+ \delta  ||u_0|| \lambda (t) E(t) +   
\\
\\  
&+ \delta  \displaystyle  \bigg[  ||y_0||+  ||w_0|| + \bigg|\frac{c}{\beta}- \frac{h}{d} \bigg| \bigg] t \, \lambda (t) E(t)+
\\
\\  
&   \displaystyle  \delta\bigg[ ||\varphi ||
 + \bigg|\frac{h}{d}- \frac{c}{\beta} \bigg|\bigg] M +(\delta d-\varepsilon \beta)\bigg( ||w_0 ||+  \frac{c}{\beta }    \bigg)  h(t)
    \end{array}
  \end{equation}                                                     
 \end{theorem}

\vspace{3mm}\noindent where constants $\,\, q,\,S, \,M,\,N $ are introduced in $(\ref{e46})_2$, (\ref{2s}), (\ref{E}), and  (\ref{EE}), respectively.

Besides, functions $ \,\, C(t),\, \,\lambda (t),\,E(t),\,L(t), \,g(t),\,h(t) $ are defined in $( \ref{a218})_3 $ (\ref{e47}),\,\,$(\ref{hEE})_{1,2},$\, (\ref{g}),\,and \,(\ref{h}).

\vspace{1mm}  According to   (\ref{abcc})  and  by means of  inequalities (\ref{439}), (\ref{a441}),  (\ref{b38}), and (\ref{bbb38}),  estimate  (\ref{68}) follows.

 \vspace{1mm} \noindent As for inequalities (\ref{fff}) and (\ref{ff}), functions defined in  (\ref{17}) have to be considered. 

\vspace{2mm} \noindent More precisely, from $(\ref{17})_1 $ and (\ref{yabcc}), taking into account inequalities (\ref{b38}), (\ref{bb38}), (\ref{bbb38}),    (\ref{ddccc}) and (\ref{ccc}), estimate (\ref{fff})  is  proved.

\vspace{1mm} \noindent Analogously, from $(\ref{17})_2 $ and (\ref{yabccc}), for (\ref{b38})-(\ref{dccc}) and (\ref{bccc}), also  (\ref{ff})  holds.

\hbox{}\hfill\rule{1.85mm}{2.82mm}


\vspace{2mm}\textbf{Remark} These estimates show  that  the solution of  the  FitzHugh -Rinzel system is bounded for all t. Besides,  when  t tends to infinity, the effect of the  non linear term $ \varphi (x,t) $ is bounded,  while   the effects of  initial perturbances $ u_0(x),  w_0(x), y_0(x) $ are vanishing.

\vspace{3mm}

\textbf{
Acknowledgements}
 
The present work has been  developed with the economic support of MIUR (Italian Ministry of University and Research) performing the  activities of the project   
ARS$01_{-}  00861$  “Integrated collaborative 
 systems for smart factory - ICOSAF".

 The  paper has been performed under the auspices of G.N.F.M. of INdAM.

  The author is grateful to the anonymous reviewers for their comments and suggestions.
 

The author declares that she has no conflict of interest.

\end{document}